
\magnification1200
\input amstex.tex
\documentstyle{amsppt}

\hsize=13cm
\vsize=18cm
\hoffset=1cm
\voffset=2cm

\footline={\hss{\vbox to 2cm{\vfil\hbox{\rm\folio}}}\hss}
\nopagenumbers
\def\DJ{\leavevmode\setbox0=\hbox{D}\kern0pt\rlap
{\kern.04em\raise.188\ht0\hbox{-}}D}

\def\txt#1{{\textstyle{#1}}}
\baselineskip=13pt
\def\hf{{\textstyle{1\over2}}}
\def\a{\alpha}\def\b{\beta}
\def\d{{\,\roman d}}
\def\e{\varepsilon}
\def\f{\varphi}
\def\G{\Gamma}

\def\s{\sigma}
\def\t{\theta}
\def\={\;=\;}

\def\D{\Delta}

\def\R{\Re{\roman e}\,} 
\def\z{\zeta}

 \def\t{\theta}
\def\hf{{\textstyle{1\over2}}}
\def\txt#1{{\textstyle{#1}}}
\def\f{\varphi}

\font\tenmsb=msbm10
\font\sevenmsb=msbm7
\font\fivemsb=msbm5
\newfam\msbfam
\textfont\msbfam=\tenmsb
\scriptfont\msbfam=\sevenmsb
\scriptscriptfont\msbfam=\fivemsb
\def\Bbb#1{{\fam\msbfam #1}}

\def \NN {\Bbb N}

\font\ff=cmr8
\font\jj=cmcsc8
\def\txt#1{{\textstyle{#1}}}
\baselineskip=13pt

\font\teneufm=eufm10
\font\seveneufm=eufm7
\font\fiveeufm=eufm5
\newfam\eufmfam
\textfont\eufmfam=\teneufm
\scriptfont\eufmfam=\seveneufm
\scriptscriptfont\eufmfam=\fiveeufm
\def\mathfrak#1{{\fam\eufmfam\relax#1}}

\font\tenmsb=msbm10
\font\sevenmsb=msbm7
\font\fivemsb=msbm5
\newfam\msbfam
     \textfont\msbfam=\tenmsb
      \scriptfont\msbfam=\sevenmsb
      \scriptscriptfont\msbfam=\fivemsb
\def\Bbb#1{{\fam\msbfam #1}}

\def \NN {\Bbb N}

  \def\rightheadline{{\hfil{\ff
  Higher moments of the error term in the divisor problem}
\hfil\tenrm\folio}}

  \def\leftheadline{{\tenrm\folio\hfil{\ff
   A. Ivi\'c and W. Zhai}\hfil}}
  \def\emptyheadline{\hfil}
  \headline{\ifnum\pageno=1 \emptyheadline\else
  \ifodd\pageno \rightheadline \else \leftheadline\fi\fi}

\topmatter
\title
HIGHER MOMENTS OF THE ERROR TERM IN THE DIVISOR PROBLEM
\endtitle
\author   Aleksandar Ivi\'c and Wenguang Zhai${}^*$ \endauthor
\address
Aleksandar Ivi\'c, Katedra Matematike,  Universitet u Beogradu,
Rudarsko-geol-\break o\v ski fakultet, \DJ u\v sina 7, 11000 Beograd, Serbia.
\medskip
Wenguang Zhai, Department of Mathematics, Shandong Normal
University, Jinan 250014, P.R. China
\bigskip
\endaddress
\keywords
The Dirichlet divisor problem, moments, mean fourth power, short intervals
\medskip
${}^*$W. Zhai is supported by the Nat. Natural
Science Foundation of China (Grant No. 10771127)
\endkeywords
\email {\tt
ivic\@rgf.bg.ac.rs, aivic\_2000\@yahoo.com
\medskip\hskip22mm
zhaiwg\@hotmail.com}
\endemail
\dedicatory
Dedicated to the memory of A.A. Karatsuba
\enddedicatory
\abstract
{It is proved that, if $k\geqslant2$ is a fixed integer and $1 \ll H \leqslant \hf X$, then
$$
\int_{X-H}^{X+H}\D^4_k(x)\d x \ll_\e X^\e\Bigl(HX^{(2k-2)/k} +
H^{(2k-3)/(2k+1)}X^{(8k-8)/(2k+1)}\Bigr),
$$
where $\D_k(x)$ is the error term in the general Dirichlet divisor problem.
The proof uses the Vorono{\"\i}--type formula
for $\D_k(x)$, and the bound of Robert--Sargos for the number of integers
when the difference of four $k$--th roots is small. We also investigate the size
of the error term in the asymptotic formula for the $m$-th moment of $\D_2(x)$.
}
\endabstract
\endtopmatter
\document

\head
1. Introduction and statement of results
\endhead
For a fixed $k\in\NN$, let
$$
\D_k(x) \;:=\; \sum_{n\leqslant x}d_k(n) - xP_{k-1}(\log x)\leqno(1.1)
$$
denote the error term in the (general) {\it Dirichlet divisor
problem} (sometimes also called the {\it Piltz divisor problem},
especially in the case when $k=3$). Here $d_k(n)$ denotes the number
of ways $n$ may be written as a product of $k$ factors (so that
$d_1(n) \equiv 1$ and $d_2(n) = d(n)$ is the number of divisors of
$n$), and  $P_{k-1}(z)$ is a suitable polynomial of degree $k-1$ in
$z$ (see e.g., [5, Chapter 13] and [12, Chapter 12]) for more
details). The function $\D_k(x)$ takes both positive and negative
values. It has finite jumps when $x = n\in\NN$ which can be of the
magnitude $\exp(C(k)\log n/\log\log n)$, the maximal order of
$d_k(n)$.

In particular, we have that
$$
\D_2(x) \equiv \D(x) =   \sum_{n\leqslant x}d(n) - x(\log x + 2\gamma - 1)
\leqno(1.2)
$$
represents the error term in the classical Dirichlet divisor problem
($\gamma = -\G'(1) = 0.5772\ldots\,$ is Euler's constant).
A vast literature exists on the estimation of $\D_k(x)$ and
especially on $\D(x)$ (op. cit.),
both pointwise and in various means. This concerns in particular various
mean square results concerning $\D(x)$.
In this work we shall be concerned with the higher moments,
especially the fourth moment of
$\D_k(x)$ in ``short" intervals of the form $[X-H,\,X+H]$,
which is the next ``natural'' moment after the square. Here ``short"
means that the relevant range for $H$ is  $H = o(X)$ as $\,X\to\infty$.

\medskip
We begin by noting that the first author in [4] (see also [5, Chapter 13]) proved
a large values estimate for $\D(x)$, which yielded the bound
$$
\int_1^X\D^4(x)\d x \;\ll_\e\; X^{2+\e},\leqno(1.3)
$$
where here and later $\e$ denotes arbitrarily small, positive
constants, which are not necessarily the same ones at each
occurrence. The asymptotic formula for the fourth moment
with an error term was obtained by K.-M. Tsang [14]. He has sharpened (1.3) to
$$
\int_1^X\D^4(x)\d x \= CX^2 + O_\e(X^{\gamma+\e})  \leqno(1.4)
$$
with an explicitly given $C\,(>0)$ and $\gamma = 45/23 = 1.956\ldots\;$.
Tsang also proved an asymptotic formula for the integral of the cube of $\D(x)$,
namely
$$
\int_1^X\D^3(x)\d x \= BX^{7/4} + O_\e(X^{\b+\e})  \leqno(1.5)
$$
with explicit $B>0$ and $\b = {47\over28} = 1.6785\ldots\;$. Later
Ivi\'c--Sargos [7] obtained the better values $\b = {7\over5} =
1.4,\, \gamma = {23\over12} = 1.91666\ldots\,$, and the second
author [15] further reduced the value of $\gamma$ to ${53\over28}= 1.8926\dots\;$.
Higher power moments of $\D(x)$ were studied by D.R. Heath--Brown [2] and
the second author. In [15] he proved that
$$
\int_1^X\D^m(x)\d x = C_mX^{1+m/4} + O_\e\Bigl(X^{1+m/4-\eta_m+\e}\Bigr)
\leqno(1.6)
$$
holds for integers $m$ satisfying $5\leqslant m\leqslant9$ with some explicit
$\eta_m>0$. He gives in closed (although complicated) form the
constants $C_m$, and it is conjectured that (1.6) holds $\forall m\in\NN$.

\medskip
Concerning the true order of $\D_k(x)$,
a classical conjecture states that $\D_k(x) \ll_{k,\e} x^{(k-1)/(2k)+\e}$, while
on the other hand $\D_k(x) = \Omega(x^{(k-1)/(2k)})$ (see e.g., [5] or [12]).
For $k=2,3$ this follows heuristically from (2.4)
and (when $k=2$) from the asymptotic formula (1.4).
The sharpest known omega-result  result for $\D_k(x)$ is due to K. Soundararajan [10],
who proved that
$$\Delta_k(x)= \Omega\Bigl\{G_k(x)\Bigr\}
\qquad(k\geqslant2)\leqno(1.7)
$$
with
$$
\eqalign{
G_k(x) :&= (x\log x)^{(k-1)/(2k)}
(\log\log x)^a (\log\log\log x)^{-b}\cr
a&= \frac{k+1}{2k} (k^{(2k)/(k+1)}-1), \quad b= \frac{3k-1}{4k}.\cr}
\leqno(1.8)
$$
Here, as usual, $f(x) = \Omega(g(x))\;(g(x) >0$ for $x\geqslant C$)
 means that $\lim_{x\to\infty}\frac{f(x)}{g(x)}=0$
does not hold. Thus there exist a sequence $\{x_n\}$ tending to infinity such
that, for some $\delta>0$,
$$
|f(x_n)| \;>\;\delta g(x_n)\qquad(n\geqslant n_0(\delta)).\leqno(1.9)
$$

\medskip
Our first aim is to prove an omega-result for the error term
in (1.6). This is contained in

\medskip
THEOREM 1. {\it If  $\eta_m$ is defined by} (1.6) {\it then, for all $m\geqslant2$,
$\eta_{m} \;\le \;3/4.$ More precisely, we have
$$
\int_1^X\D^m(x)\d x = C_mX^{1+m/4} + \Omega\Bigl(G^{m+1}(X)\log^{-1}X\Bigr),
\leqno(1.10)
$$
where $G(x) \equiv G_2(X)$ is given by} (1.8).
\medskip
{\bf Remark 1}. It would be interesting to investigate the sign of $C_m$ in
(1.6) in the general case; the values for $2\leqslant m \leqslant 9$
are all positive.

{\bf Remark 2}. A result analogous to Theorem 1 could be obtained for
the $m$-th moment of the general error-term function
$\D_k(x)$. However, except for the asymptotic formula
$$
\int_1^X\D_3^2(x)\d x = (10\pi^2)^{-1}\sum_{n=1}^\infty d_3^2(n)n^{-4/3}
X^{5/3} + O_\e(X^{14/9+\e})
$$
of K.-C. Tong [13], there are no other asymptoptic
formulas for moments of $\D_k(x)$
when $k>2$. Hence such a result at present would not have much practical value.

\smallskip
Moments of $\D(x)$ over short intervals were
investigated by Lau--Tsang [8]. In particular, they give the existence
of
$$
\lim_{X\to\infty}{1\over HX^{m/4}}\int_{X-H}^{X+H}|\D(x)|^m\d x\qquad(m\in\NN)
$$
under certain conditions on $m$ and $H$.
Lau and Tsang also investigated the above
integral with $\D^m(x)$ (i.e., without the absolute values).
For $k>2$ there seem to be no analogous results available for $\D_k(x)$.

\smallskip
Our result, which is primarily of significance when $2 \leqslant k \leqslant 4$,
deals with the fourth moment of $\D_k(x)$ in short intervals. It is  the following

\medskip
THEOREM 2. {\it If $k\geqslant2$ is a fixed integer and $1 \ll H \leqslant \hf X$, then}
$$
\int_{X-H}^{X+H}\D^4_k(x)\d x \;\ll_\e\; X^\e\Bigl(HX^{(2k-2)/k} +
H^{(2k-3)/(2k+1)}X^{(8k-8)/(2k+1)}\Bigr). \leqno(1.11)
$$

\medskip
{\bf Remark 3}. Note that, for $H\leqslant\hf X$ and $k\geqslant 3$ we
have
$$
HX^{(2k-2)/k} \;\leqslant\; H^{(2k-3)/(2k+1)}X^{(8k-8)/(2k+1)},
$$
and for $k = 2$ this holds when $H\leq X^{3/4}$.
The importance of $HX^{(2k-2)/k+\e}$ is that it is the ``expected" order
of the integral in (1.11), in view of the conjecture
$\D(x) \ll_{k,\e} x^{(k-1)/(2k)+\e}$.
\medskip
{\bf Remark 4}. One of the reasons that we treat
the fourth power of  $\D_k(x)$ is the
result of Robert--Sargos [9], embodied in Lemma 2. The case of the odd
moments is more difficult, since $\D_k^m(x)$ takes both positive
and negative values if $m$ is odd.

\medskip
{\bf Remark 5}. In the case $k=2$ the Theorem 2 yields
$$
\int_{X-H}^{X+H}\D^4(x)\d x \;\ll_\e\;
 X^\e\Bigl(HX + X^{8/5}H^{1/5}\Bigr)\qquad(1\ll H\leqslant\hf X),
\leqno(1.12)
$$
while (1.4) yields
$$
\int_{X-H}^{X+H}\D^4(x)\d x \;\ll_\e\; HX + X^{53/28+\e}.
\leqno(1.13)
$$
Note that (1.12) improves (1.13). However, in this case one can
obtain an even sharper bound. Let
$$
\a \;:=\; \inf\Bigl\{\;a \;:\; \D(x) \ll x^a\;\Bigr\}.\leqno(1.14)
$$
Then
$$
\int_{X-H}^{X+H}\D^4(x)\d x \;\ll_\e\; X^\e(HX + X^{1+2\a})
\qquad(\sqrt{X} \ll H \ll X).\leqno(1.15)
$$
The bound (1.15) follows from Theorem 6 of Lau--Tsang [8] when
$A=4$. It is proved by employing the large value method of [4] and
[5, Chapter 13], in particular see (13.52) of [5].

\medskip
{\bf Remark 6}.  It is well-known (i.e., follows from (1.4)) that
$\a\geqslant1/4$, hence the bound in (1.15) trivially holds when $H
\ll \sqrt{X}$ . On the other hand the best upper bound for $\a$ in
(1.14) at present is
$$
\a \; \leqslant\; {131\over416} = 0.314903\ldots\,,\leqno(1.16)
$$
due to M.N. Huxley [3]. Hence combining (1.15) and (1.16) we obtain
in fact
$$
\int_{X-H}^{X+H}\D^4(x)\d x \;\ll_\e\; X^\e(HX + X^{339/208})\qquad(1 \ll H \ll X).
$$

\head
2. The necessary Lemmas
\endhead

 We begin with  the elementary

\medskip
LEMMA 1. {\it For $1 \ll H \ll X$ and any integer $k\geqslant2$ we have}
$$
\D_k(X) = {1\over H}\int_{X}^{X+H}\D_k(x)\d x + O_{\e,k}(HX^\e).
\leqno(2.1)
$$

\medskip
{\bf Proof}. Since $d_k(n) \ll_{\e,k} n^\e$, it follows from the
defining relation (1.1) that
$$\eqalign{\cr&
\D_k(X) - {1\over H}\int_{X}^{X+H}\D_k(x)\d x\cr& = {1\over
H}\int_{X}^{X+H}(\D_k(X) - \D_k(x))\d x\cr& \ll {1\over
H}\int_{X}^{X+H}\Bigl(\bigl|\sum_{X\leqslant n\leqslant
x}d_k(n)\bigr| + O(HX^\e)\Bigr)\d x \cr& \ll_\e HX^\e,\cr}
$$
which gives (2.1). By using the well-known result of P. Shiu [11] on
multiplicative functions, one can improve the error term in (2.1)
to $O_k(H\log^{k-1}X)$ in the range $X^\e \leqslant H \ll X$.
\medskip
LEMMA 2. {\it Let $k\geqslant 2$ be a fixed
integer and $\delta > 0$ be given.
Then the number of integers $\;n_1,n_2,n_3,n_4\;$ such that
$\;3\leqslant N < n_1,n_2,n_3,n_4 \leqslant 2N\;$ and}
$$
\left|n_1^{1/k} + n_2^{1/k} - n_3^{1/k} - n_4^{1/k}\right| < \delta N^{1/k}
$$
{\it is, for any given $\e>0$,}
$$
\ll_\e\; N^\e(N^4\delta + N^2).\leqno(2.2)
$$

\medskip
This result was proved by O. Robert--P. Sargos [9]. It represents an
arithmetic tool which is useful in dealing with various analytic
problems.

It seems reasonable to conjecture that, for fixed $\ell\geqslant2$, the number
of integers $n_1,\ldots,n_\ell,n_{\ell+1},\ldots\,, n_{2\ell}$ for which
$$
\left|n_1^{1/k} +\ldots +n_\ell^{1/k} - n_{\ell+1}^{1/k} -\ldots- n_{2\ell}^{1/k}\right|
< \delta N^{1/k}
$$
holds is, for any given $\e>0$,
$$
\ll_\e\; N^\e(N^{2\ell}\delta + N^\ell).\leqno(2.3)
$$
This conjecture is very strong, and already the truth of (2.3) for $\ell=3$ would
allow one to treat the sixth power of $\D_k(x)$ in short intervals.

\medskip
LEMMA 3. {\it For fixed $k\geqslant2$ and $1\ll N \ll x$, we have}
$$
\eqalign{
\D_k(x) & = {x^{(k-1)/(2k)}\over\pi\sqrt{k}}\sum_{n\leqslant N}
d_k(n)n^{-(k+1)/(2k)}\cos\Bigl(2k\pi(xn)^{1/k}+ \txt{1\over4}(k-3)\pi\Bigr)\cr&
+ O_\e\Bigl(x^{(k-1)/k+\e}N^{-1/k}\Bigr).
\cr}\leqno(2.4)
$$

\medskip {\bf Proof}.   The explicit, Vorono{\"\i}--type formula (2.4), is well
known in the case when $k=2$ (see [5, Chapter 3] or [12, Chapter 12]) for a proof).
However, in the general case it does not seem to appear in the literature
and a proof (based on the classical proof in the case when $k=2$) is in order.

To begin with note that,
for $x^\e\leqslant T\ll x$, Perron's inversion formula (see e.g., the Appendix of [5])
gives
$$
\sum_{n\leqslant x}d_k(n) = {1\over2\pi i}\int_{1+\e-iT}^{1+\e+iT}
\z^k(s)x^ss^{-1}\d s + O_\e(x^{1+\e}T^{-1}).
$$
We replace the segment of integration  by the segment $[-\e - iT,\,-\e + iT]$, passing
over the pole of $\z^k(s)$ at $s=1$. By the residue theorem this yields the
term $xP_{k-1}(\log x)$ (cf. (1.1)). The horizontal segments
$[-\e\pm iT,\,1+\e\pm iT]$ make a contribution which is
$$
\ll \int_{-\e}^{1+\e}x^\s T^{-1}|\z(\s+iT)|^k\d \s \ll_\e x^\e(xT^{-1} + T^{k/2-1}).
$$
This is obtained by using the standard convexity bound (see [5, Chapter 1])
$$
\z(\s+it) \ll t^{(1-\s)/3}\log t \qquad(\hf \leqslant \s\leqslant1,\;t \gg 1)
$$
and the functional equation for $\z(s)$ in the form
$$
\z(s) = \chi(s)\z(1-s),\; \chi(s) = \Bigl({2\pi\over t}\Bigr)^{\s+it-1/2}{\roman e}
^{i(t+\pi/4)}\Bigl(1+ O\Bigl({1\over t}\Bigr)\Bigr)\quad(t\geqslant t_0 > 0).
$$
It follows that
$$
\D_k(x) = 2\R I + O_\e\Bigl(x^\e(xT^{-1} + T^{k/2-1})\Bigr),\leqno(2.5)
$$
where we have set
$$
I := {1\over2\pi}\int_1^T\chi^k(-\e+it)\sum_{n=1}^\infty d_k(n)n^{-1-\e+it}
x^{-\e+it}{\d t\over-\e+it},\leqno(2.6)
$$
and the series in (2.6) is absolutely convergent. Therefore we may change the order
of integration and summation to obtain
$$
\eqalign{
I&= {1\over2\pi}\int_1^T\chi^k(-\e+it)\sum_{n=1}^\infty d_k(n)n^{-1-\e+it}
x^{-\e+it}{\d t\over it} + O_\e(x^\e T^{k/2-1})\cr&
= {{\roman e}^{k\pi i/4}\over 2\pi ix^\e}\sum_{n=1}^\infty d_k(n)n^{-1-\e}
\int_1^T\Bigl({2\pi\over t}\Bigr)^{-k\e+kit-k/2}{\roman e}^{kit}(xn)^{it}
{\d t\over t} + O_\e(x^\e T^{k/2-1})\cr&
= {{\roman e}^{k\pi i/4}(2\pi)^{-k\e-k/2}\over 2\pi ix^\e}\sum_{n=1}^\infty {d_k(n)\over
n^{1+\e}}
\int\limits_1^T t^{k\e+k/2-1}{\roman e}(F(t))\d t +
O_\e(x^\e T^{k/2-1})\cr}\leqno(2.7)
$$
with ${\roman e}(z) = {\roman e}^{2\pi iz}$, and
$$
2\pi F(t) = 2\pi F(t;x,k,n) := -kt\log{t\over2\pi} + kt + t\log(xn).
$$
Thus the saddle point of the last exponential integral in (2.7) (root of $F'(t) = 0$)
is
$
t = t_0 = 2\pi(xn)^{1/k},
$
so that $t_0 \leqslant T$ for $n\leqslant (T/(2\pi))^kx^{-1}$. Therefore
the parameter $T$ is determined to satisfy
$$
N + {1\over2} \;=\; \left({T\over2\pi}\right)^kx^{-1}\qquad(N\in\NN),
$$
and $T\ll x$ holds because $N\ll x, k\geqslant2$.
The exponential integral in (2.7) is then evaluated e.g., by the well-known result of F.V.
Atkinson [1] (see also [5, Chapter 2] or [6, Chapter 2]). The main term
for the integral in the last line in the expression (2.7)
for $I$ (see p. 65 of [5]) is, since $|F''(t_0)|^{-1/2} = 2\pi(xn)^{1/(2k)}k^{-1/2}$,
$F(t_0) = k(xn)^{1/k}$,
$$
(2\pi)^{k\e+{k\over2}-1}k^{-{1\over2}}(xn)^{\e+{1\over2}-{1\over k}}2\pi(xn)^{1/(2k)}
{\roman e}( k(xn)^{1/k}){\roman e}^{-{\pi i\over4}}.
$$
For $n>N, 1\leqslant t \leqslant T$ we have
$$
F'(t) \;= \;{1\over2\pi}\log\Biggl({xn\over\Bigl({t\over2\pi}\Bigr)^k}\Biggr)
\;\geqslant\; {1\over2\pi}\log {n\over N+\hf},
$$
so that by the first derivative test (i.e., Lemma 2.1 of [5]) the contribution of
$n>N$ to $I$ is
$$
\eqalign{
&\ll_\e x^{-\e}T^{k/2-1}\Biggl(\sum_{N<n\leqslant2N}d_k(n)n^{-1-\e}
\cdot{1\over\log{n\over N+\hf}}
+ \sum_{n>2N}d_k(n)n^{-1-\e}\Biggr) + x^\e T^{k/2-1}\cr&
\ll_\e x^\e T^{k/2-1}\cr}
\leqno(2.8)
$$
on writing $n = N+r, 1\leqslant r \le N$ in the first sum above.
There remain the error terms in Atkinson's formula (cf. (2.16) of [5]) for $n\leqslant N$.
The first two error terms are clearly absorbed by the error terms in (2.5) and so
is the third (corresponding to $a=1$). The fourth  error term makes a contribution
which is
$$
\eqalign{&
\ll_\e x^\e T^{k/2-1}\sum_{n\leqslant N}d_k(n)n^{-1-\e}\cdot{1\over\log\Bigl({xn\over
\bigl({T\over2\pi}\bigr)^k}\Bigr)}\cr&
\ll_\e x^\e T^{k/2-1}\Bigl(\sum_{n\leqslant N/2} + \sum_{N/2\leqslant n \leqslant N}\Bigr)
\ll_\e x^\e T^{k/2-1}\cr}
$$
on proceeding similarly as in (2.8). Since $T\asymp (xN)^{1/k}$, the assertion
of the lemma follows.

\head
3.  Proof of  Theorem 1
\endhead

We shall use the following version of Lemma 1: 
$$
\D(x) = {1\over H}\int_X^{X+H}\D(x)\d x + O(H\log X)\qquad(X^\e \leqslant H \leqslant X),
 \leqno(3.1)
$$
which follows, as was mentioned in connection with Lemma 1, from  P. Shiu's   bound (see [11])
$$
\sum_{x<n\le x+h}d(n) \;\ll\; h\log x\qquad(x^\e \le h \le x).\leqno(3.2)
$$
Namely as in the proof of Lemma 1 we have, for $X^\e \leqslant H \leqslant X$,
$$\eqalign{\cr&
\left|\D(X) - {1\over H}\int_{X}^{X+H}\D(x)\d x\right| = \left|{1\over
H}\int_{X}^{X+H}(\D(X) - \D(x))\d x\right|\cr& \le {1\over
H}\int_{X}^{X+H}\Bigl(\sum_{X\leqslant n\leqslant
x}d(n)+ O(H\log X)\Bigr)\d x \cr&
\le {1\over H}\int_{X}^{X+H}\Bigl(\sum_{X\leqslant n\leqslant
X+H}d(n) + O(H\log X)\Bigr)\d x \ll H\log X,
\cr}
$$
where we used (1.1) with $k=2$ and (3.2).

\smallskip
 Now we take $X = x_n$, the point for which the omega-result
(1.8)--(1.9) is attained when $k=2$. As Soundararajan's method of
proof does not tell whether $\D(x_n)$ is positive or negative, we
choose $\t = {\roman {sgn}}\, X = 1$ if $\D(X)>0$, and $\t = -1$ it
$\D(X)<0$. Hence $\t \D(X) = |\D(X)|.$
When $X+X^\e\leq x\leq X+H$ we have $|X-x| \geq X^\e$, hence
$|\D(X) - \D(x)| \leqslant B_1|X-x|\log X$ by (3.2),
similarly as in the proof (3.1).
If $X\leq x \leq X+X^\e$  we use the the trivial estimate $d(n)\ll_\e n^\e$
and (1.1). Hence we obtain
$$
|\D(X) - \D(x)| \leqslant B_1|X-x|\log X+X^\e\qquad(X\leq x \leq X+H)\leqno(3.3)
$$
for some $B_1>0$. Therefore if we take
$$
H \;=\; B_2\,{G(X)\over\log X}\leqno(3.4)
$$
with a   sufficiently small constant $B_2>0$, it follows from (3.3)
that $\t\D(x) >0$ for $X  \leqslant x \leqslant X+H$. On multiplying
by $\t$ we obtain  from (3.1), for some $B_3>0$,
$$
B_3(X\log X)^{1\over4}(\log_2 X)^{{3\over4}(2^{4/3}-1)}(\log_3
X)^{-{5\over8}} = B_3G(X) \leqslant {1\over H}\int_X^{X+H}\t\D(x)\d
x,\leqno(3.5)
$$
where $\log_r X = \log(\log_{r-1}X)$. We raise (3.5) to the $m$-th power,
and use H\"older's inequality for integrals, since the integrand in (3.5) is positive.
It follows that
$$
\eqalign{ \Bigl(B_3G(X)\Bigr)^m & \leqslant {1\over
H^m}\Bigl(\int_{X }^{X+H}\t\D(x)\d x \Bigr)^m \leqslant {1\over
H}\int_{X }^{X+H}\t^m \D^m(x)\d x\cr& = {1\over H}\Bigl(\t^m E_m(X+H) -
\t^m E_m(X ) + O(HX^{m/4})\Bigr),\cr} \leqno(3.6)
$$
where (cf. (1.6))
$$
E_m(X) \;:=\; \int_1^X \D^m(x)\d x - C_m X^{1+m/4}.
$$
Should it happen that, for sufficiently small $c_0 >0$ and $X \geqslant X_0$,
$$
|E_m(X)| \;\leqslant\; c_0{G^{m+1}(X)\over\log X},
$$
then it follows from (3.4) and (3.6) that, for sufficiently large $X$,
$$
\hf \Bigl(B_3G(X)\Bigr)^m \leqslant 2B_2^{-1}c_0G^m(X),
$$
which is a contradiction if $c_0 < {1\over4}B_2B_3^{m}$. This proves Theorem 1,
and $\eta_m \leqslant 3/4$ is a consequence of the explicit expression
for $G(x)$ in (3.5).

\head
3.  Proof of  Theorem 2
\endhead
\medskip
We pass now to the proof of Theorem 2.
Let $\f(x)\;(\geqslant0)$ be a smooth function supported in $[X-2H,X+2H]$ such
that $\f(x) = 1$ for $x\in[X-H,X+H]$, so that
$$
\f^{(r)}(x) \;\ll_r\; H^{-r}\qquad(r = 0,1,2,\ldots).
\leqno(4.1)
$$
If $x\asymp X,N\ll X$, then from Lemma 3 we obtain
$$
\eqalign{&
\int_{X-H}^{X+H}\D^4_k(x)\d x \leqslant \int_{X-2H}^{X+2H}\f(x)\D^4_k(x)\d x\cr&
\ll_\e X^{(2k-2)/k+\e}\max_{K\leqslant N}\int_{X-2H}^{X+2H}\f(x)
\Bigl|\sum_{K<n\leqslant K'\leqslant2K}d_k(n)n^{-(k+1)/(2k)}{\roman e}^{2k\pi i(xn)^{1/k}}
\Bigr|^4\d x\cr&
+ HX^{(4k-4)/k+\e}N^{-4/k}.\cr}\leqno(4.2)
$$
The integral in (4.2) is equal to
$$
\sum_{K<m,n,j,\ell\leqslant K'}d_k(m)d_k(n)d_k(j)d_k(\ell)(mnj\ell)^{-(k+1)/(2k)}
\int_{X-2H}^{X+2H}\f(x){\roman e}^{iDx^{1/k}}\d x,\leqno(4.3)
$$
where $(m,n,j,\ell) \in \NN^4$ and
$$
D = D_k(m,n,j,\ell) := 2k\pi\Bigl(m^{1/k} + n^{1/k} - j^{1/k} - \ell^{1/k}\Bigr).
$$
Integration by parts shows that the integral in (4.3) is
$$
{ik\over D}\int_{X-2H}^{X+2H}\Bigl(\f'(x)x^{1-1/k}+ (1 - 1/k)\f(x)x^{-1/k}\Bigr)
{\roman e}^{iDx^{1/k}}\d x.
$$
This shows that we have obtained the same type of exponential integral,
only the integrand is smaller by a factor of
$$
{1\over D}\left(H^{-1}X^{1-1/k} + X^{-1/k}\right)\; \ll\; {X^{1-1/k}\over HD}.
$$
Therefore, if we perform integration by parts $r = r(A,\e)$ times,
then in view of (4.1) we see that the contribution
of $D$ for which
$$
|D| \;>\; X^{1-1/k+\e}H^{-1}
$$
will be smaller than $X^{-A}$ for any given $A>0$. In the case when
$$
|D| \;\leqslant\; X^{1-1/k+\e}H^{-1}
$$
we can use Lemma 2 (with $\delta = X^{1-1/k}H^{-1}K^{-1/k}$)
and trivial estimation to infer that the expression in (4.3) is
$$
\eqalign{&
\ll_\e X^\e K^{-(2k+2)/k}H(K^4X^{1-1/k}H^{-1}K^{-1/k} + K^2)
\cr&
\ll_\e X^\e H + X^{1-1/k+\e}N^{2-3/k}.\cr}
$$
This gives, in view of (4.2),
$$
\int_{X-H}^{X+H}\D^4_k(x)\d x \ll_\e X^\e\left(HX^{(2k-2)/k}
+ X^{3-3/k}N^{2-3/k} + HX^{4-4/k}N^{-4/k}\right).
\leqno(4.4)
$$
The terms containing $N$ in (4.4) are equal if
$$
N \;=\; X^{(k-1)/(2k+1)}H^{k/(2k+1)}\quad(< X).
$$
Therefore (1.8) follows from (4.4), and the proof of the Theorem
is complete.

\medskip
{\bf Remark 7}. From Lemma 1 and H\"older's inequality for integrals
it follows that
$$
\D_k^4(X) \;\ll_\e\; {1\over H}\int_{X-H}^{X+H}\D^4_k(x)\d x + X^\e
H^4. \leqno(4.5)
$$
If we take $H = X^{(k-1)/(k+1)}$ in (4.5) and apply (1.8) of  the Theorem,
we obtain the known bound (follows also from Lemma 3)
$$
\D_k(X) \;\ll_\e\; X^{(k-1)/(k+1)+\e}.
$$

\bigskip

\vfill
\eject
\topglue1cm
\vskip1cm
\Refs
\bigskip

\item{[1]} {\jj F.V. Atkinson}, `The mean value of the Riemann zeta-function',
{\it Acta Math.} {\bf81}(1949), 353-376.

\item{[2]} {\jj D.R. Heath-Brown}, `The distribution of moments
in the Dirichlet divisor problems', {\it Acta Arith.} {\bf60}(1992),
389-415.

\item{[3]} {\jj M.N. Huxley},
`Exponential sums and lattice points III',
Proc. London Math. Soc., (3) {\bf 87}(2003), 591-609.

\item{[4]} {\jj A. Ivi\'c}, `Large values of the error term in the divisor problem',
    {\it Inventiones Math}. {\bf71} (1983), 513-520.

\item{[5]} {\jj A. Ivi\'c}, `The Riemann zeta-function', John Wiley \&
Sons, New York, 1985 (2nd ed. Dover, Mineola, New York, 2003).

\item{[6]} {\jj A. Ivi\'c}, `The mean values of the Riemann zeta-function',
LNs {\bf 82}, Tata Inst. of Fundamental Research, Bombay (distr. by
Springer Verlag, Berlin etc.), 1991.

\item{[7]} {\jj A. Ivi\'c and P. Sargos},
`On the higher power moments of the error term in the
divisor problem', {\it Illinois J.Math.} {\bf81}(2007), 353-377.

\item{[8]} {\jj Y.-K. Lau and K.-M. Tsang}, `Moments over short
intervals', {\it Arch. Math.} {\bf84}(2005), 249-257.

\item{[9]} {\jj O. Robert and P. Sargos},  `Three-dimensional
exponential sums with monomials', {\it J. reine angew. Math.} {\bf591}(2006), 1-20.

\item{[10]} {\jj K. Soundararajan},
Omega results for the divisor and circle problems,
J. Int. Math. Res. Not. 2003, No. {\bf36}, 1987-1998(2003).

\item{[11]} {\jj P. Shiu}, `A Brun-Titchmarsh theorem for multiplicative functions',
{\it J. reine Angew. Math.} {\bf31}(1980), 161-170.

\item{[12]} {\jj E.C. Titchmarsh}, `The theory of the Riemann zeta-function'
(2nd ed.),  University Press, Oxford, 1986.

\item{[13]} {\jj K.-C. Tong}, On divisor problems III (Chinese),
Acta Math. Sinica {\bf6}(1956), 515-541.

\item{[14]} {\jj K.-M. Tsang}, `Higher power moments of $\D(x)$, $E(t)$
and $P(x)$', {\it Proc. London Math. Soc. }(3){\bf 65}(1992), 65-84.

\item{[15]} {\jj W. Zhai}, `On higher-power moments of $\D(x)$' {\it Acta Arith.}
{\bf112}(2004),
367-395,  Part II, ibid. {\bf114}(2004), 35-54, Part III, ibid. {\bf118}(2005), 263-281,
and Part IV  (in Chinese),
Acta Math. Sin., Chin. Ser. {\bf49}, No. 3, (2006), 639-646.

\endRefs

\enddocument

\bye